\documentclass[12pt,reqno]{amsart}
\usepackage{eufrak}
\usepackage{amsmath,amsthm,amsfonts,amssymb,times}
\usepackage{graphicx,pgfarrows,pgfnodes}
\usepackage[mathscr]{eucal}
\usepackage{verbatim}
\usepackage{url,cite,mleftright}
\usepackage{tkz-fct}
\usepackage{tikz}
\usepackage[
pdfauthor={},
pdfkeywords={},
pdftitle={},
pdfcreator={},
pdfproducer={},
linktocpage,colorlinks,bookmarksnumbered,linkcolor=blue,
citecolor=red,urlcolor=red]{hyperref}
\usepackage{stackrel}

\DeclareMathOperator{\ab}{ab}
\DeclareMathOperator{\Aut}{Aut}
\DeclareMathOperator{\BGL}{BGL}
\DeclareMathOperator{\Br}{Br}
\DeclareMathOperator{\card}{card}
\DeclareMathOperator{\ch}{ch}
\DeclareMathOperator{\Char}{char}
\DeclareMathOperator{\CHur}{CHur}
\DeclareMathOperator{\Cl}{Cl}
\DeclareMathOperator{\coker}{coker}
\DeclareMathOperator{\Conf}{Conf}
\DeclareMathOperator{\disc}{disc}
\DeclareMathOperator{\End}{End}
\DeclareMathOperator{\et}{\text{\'et}}
\DeclareMathOperator{\Fix}{Fix}
\DeclareMathOperator{\Gal}{Gal}
\DeclareMathOperator{\GL}{GL}
\DeclareMathOperator{\Hom}{Hom}
\DeclareMathOperator{\Hur}{Hur}
\DeclareMathOperator{\im}{im}
\DeclareMathOperator{\Ind}{Ind}
\DeclareMathOperator{\Inn}{Inn}
\DeclareMathOperator{\Irr}{Irr}
\DeclareMathOperator{\lcm}{lcm}
\DeclareMathOperator{\Mor}{Mor}
\DeclareMathOperator{\ord}{ord}
\DeclareMathOperator{\Out}{Out}
\DeclareMathOperator{\Perm}{Perm}
\DeclareMathOperator{\PGL}{PGL}
\DeclareMathOperator{\Pin}{Pin}
\DeclareMathOperator{\PSL}{PSL}
\DeclareMathOperator{\rad}{rad}
\DeclareMathOperator{\SL}{SL}
\DeclareMathOperator{\SO}{SO}
\DeclareMathOperator{\Li}{Li}
\DeclareMathOperator{\Spec}{Spec}
\DeclareMathOperator{\Spin}{Spin}
\DeclareMathOperator{\St}{St}
\DeclareMathOperator{\Surj}{Surj}
\DeclareMathOperator{\Syl}{Syl}
\DeclareMathOperator{\tame}{tame}
\DeclareMathOperator{\Tr}{Tr}

\newcommand{\eps}{\varepsilon}
\newcommand{\QED}{\hspace{\stretch{1}} $\blacksquare$}
\renewcommand{\AA}{\mathbb{A}}
\newcommand{\CC}{\mathbb{C}}
\newcommand{\EE}{\mathbb{E}}
\newcommand{\FF}{\mathbb{F}}
\newcommand{\HH}{\mathbb{H}}
\newcommand{\NN}{\mathbb{N}}
\newcommand{\OO}{\mathbb{O}}
\newcommand{\PP}{\mathbb{P}}
\newcommand{\QQ}{\mathbb{Q}}
\newcommand{\RR}{\mathbb{R}}
\newcommand{\ZZ}{\mathbb{Z}}
\newcommand{\bfm}{\mathbf{m}}
\newcommand{\mcA}{\mathcal{A}}
\newcommand{\mcG}{\mathcal{G}}
\newcommand{\mcH}{\mathcal{H}}
\newcommand{\mcM}{\mathcal{M}}
\newcommand{\mcN}{\mathcal{N}}
\newcommand{\mcO}{\mathcal{O}}
\newcommand{\mcP}{\mathcal{P}}
\newcommand{\mcQ}{\mathcal{Q}}
\newcommand{\mfa}{\mathfrak{a}}
\newcommand{\mfb}{\mathfrak{b}}
\newcommand{\mfc}{\mathfrak{c}}
\newcommand{\mfI}{\mathfrak{I}}
\newcommand{\mfM}{\mathfrak{M}}
\newcommand{\mfm}{\mathfrak{m}}
\newcommand{\mfo}{\mathfrak{o}}
\newcommand{\mfO}{\mathfrak{O}}
\newcommand{\mfP}{\mathfrak{P}}
\newcommand{\mfp}{\mathfrak{p}}
\newcommand{\mfq}{\mathfrak{q}}
\newcommand{\mfz}{\mathfrak{z}}
\newcommand{\msC}{\mathscr{C}}
\newcommand{\msP}{\mathscr{P}}
\newcommand{\AGL}{\mathbb{A}\GL}
\newcommand{\Qbar}{\overline{\QQ}}
\renewcommand{\qedsymbol}{$\blacksquare$}
\newcommand{\bbone}{\mathbbm{1}}
\newcommand{\B}{B}

\DeclareMathOperator{\Td}{Todd}
\DeclareMathOperator{\Tdm}{\widetilde{Todd}}
\DeclareMathOperator{\Dop}{D}
\newcommand{\BN}{B}
\newcommand{\Note}[1]{\textcolor{red}{#1}}

\newcommand{\thre}{{\sqrt{3}}}
\newcommand{\qq}{{\(\dfrac q{\qu}\)^2}}
\renewcommand{\(}{\left\(}
\renewcommand{\)}{\right\)}
\renewcommand{\pmod}[1]{\,(\textup{mod}\,#1)}

\theoremstyle{plain}
\newtheorem{thm}{Theorem}

\newtheorem{cor}[thm]{Corollary}
\newtheorem{prop}[thm]{Proposition}

\theoremstyle{definition}
\newtheorem{defn}[thm]{Definition}

\theoremstyle{remark}

\numberwithin{equation}{section}
\numberwithin{thm}{section}
\allowdisplaybreaks
\DeclareMathOperator{\Imp}{Im}

\newcommand{\tw}[1]{\textcolor{red}{[Tanay: #1]}}
\newcommand{\p}[1]{\textcolor{teal}{[Parth: #1]}}

\begin{document}
	\title{Hurwitz Zeta Functions and Ramanujan's Identity for Odd Zeta Values}
	\author{Parth Chavan}
	\begin{abstract}
		Inspired by a famous formula of Ramanujan for odd zeta values, we prove an analogous formula involving the Hurwitz zeta function. We introduce a new integral kernel related to the Hurwitz zeta function, generalizing the integral kernel associated to Ramanujan's identity. We also derive several infinite families of identities analogous to Ramanujan's formula.

	\end{abstract}
	\subjclass[2020]{Primary $11$Mxx, $11$B$68$ and $11$F$03$}
	
	\maketitle
	\section{Introduction}
	
	The Riemann zeta function $\zeta(s)$, defined by the series
	$$\zeta(s) = \sum_{n=1}^{\infty}\frac{1}{n^s},$$
	with $\mathfrak{R}(s)>1$, is one of the most important special functions of mathematics. The critical strip $0 < \mathfrak{R}(s) < 1$ is undoubtedly the most famous region in the complex plane on account of the unsolved problem regarding the location of non-trivial zeros of $\zeta(s)$, namely, the Riemann Hypothesis, however the right-half plane $\mathfrak{R}(s) >1$ also has its own share of interesting unsolved problems, such as the arithmetic nature of odd zeta values. 
	
	One of the identities given by Ramanujan that has attracted the attention of many mathematicians over the years is the following intriguing identity involving the
	odd values of the Riemann zeta function:
	
	\begin{thm}[Ramanujan's formula for $\zeta(2n+1)$]\label{thm}
		If $\alpha$ and $\beta$ are positive real numbers such that $\alpha\beta=\pi^{2}$
		and if $n \in \mathbb{Z}\setminus{\{0\}}$, then we have 
		\begin{align}
			\alpha^{-n}\left\{ \dfrac{1}{2}\,\zeta(2n+1)+\sum_{m=1}^{\infty}\dfrac{m^{-2n-1}}{e^{2\alpha m}-1}\right\} &-\left(-\beta\right)^{-n}\left\{ \dfrac{1}{2}\,\zeta(2n+1)+\sum_{m=1}^{\infty}\dfrac{m^{-2n-1}}{e^{2\beta m}-1}\right\}\nonumber 
			\\&=2^{2n}\sum_{k=0}^{n+1}\dfrac{\left(-1\right)^{k-1}B_{2k}\,B_{2n-2k+2}}{\left(2k\right)!\left(2n-2k+2\right)!}\,\alpha^{n-k+1}\beta^{k}.\label{eq:Ramanujan_main}
		\end{align}
		where $B_{n}$ denotes the $n$-th Bernoulli number. 
	\end{thm}
	
	Theorem \ref{thm} appears as Entry 21 in Chapter 14 of Ramanujan's second
	notebook \cite{Bruce}. It also appears in a formerly unpublished
	manuscript of Ramanujan that was published in its original handwritten form with his lost notebook \cite{LostNotebook}.  For an elementary proof of Theorem \ref{thm} we refer the reader to \cite{Chavan}. For history and developments related to Ramanujan's formulas we refer the reader to \cite{Berndt}. The first published proof of Theorem \ref{thm} is due to S.L. Marulkar \cite{SL}, although he was not aware that this formula can be found in Ramanujan’s Notebooks. 
	The function $\frac{1}{e^{2\pi x}-1}$ appears in several Ramanujan's identities and has the integral representation 
	$$\frac{1}{e^{2\pi x}-1} = \int_{(c)} \frac{\zeta(1-s)}{2\cos\left(\frac{\pi s}{2}\right)}  x^{-s} \mathrm{d}s,$$
	where $(c)$ denoted the vertical line line $\Re(s)=c$ and $c$ is an arbitrary real number such that $c>1$. We call this function \textit{Ramanujan's kernel}. In this article we generalize introduce a two parameter generalization of Ramanujan's kernel and obtain several Ramanujan type identities. One reason to study this kernel is to obtain more information on the arithmetic nature of $\zeta(2m+1)$. By including a free parameter $a$, one can derive other results involving the zeta function by differentiating or integrating against $a$.
	
	Ramanujan's kernel has simple poles at $x=0$ and $x = \pm in , n \in \NN$  with residue $\frac{1}{2\pi}$ at $0,\pm in$ and thus has the partial fraction expansion 
	\[\frac{1}{e^{2\pi x}-1} = -\frac{1}{2} +\frac{1}{2\pi x} + \frac{x}{\pi} \sum_{n=1}^{\infty} \frac{1}{n^2+x^2}.\] We generalize this result to a meromorphic function with simple poles at $x=0$,$x=e^{\frac{(2j+1)i\pi }{2k}}(n+a)$ where $j\in\{0,1,\ldots,2k-1\}$, $n \in \NN$ and $a\in\CC$ with residue $\frac{1}{2k\pi }$ at $e^{\frac{(2j+1)i\pi }{2k}}(n+a)$ defined as follows
	\begin{defn}\label{kernel} Let $x \in \RR^{+} , a \in \CC$ and $k \in \NN$. Define the \textit{Hurwitz kernel} by
		\begin{align}
			\Psi\left(x,a;k\right)&:=\int_{\left(c\right)} \frac{\zeta(1-s,a)}{2k\cos\left(\frac{\pi(s+k-1)}{2k}\right)}x^{-s}\mathrm{d}s\nonumber
			\\&=\frac{2a-1}{2\pi x}-\frac{1}{2k\cos\left(\frac{\pi\left(k-1\right)}{2k}\right)} + \frac{1}{\pi }\sum_{n=0}^{\infty}\frac{x^{2k-1}}{x^{2k}+(n+a)^{2k}},
		\end{align}
	\end{defn}
	where $\zeta(s,a)$ is the Hurwitz zeta function defined in Section \ref{mainresults}. Note that $\Psi(x,1;1)$ is Ramanujan's kernel. A. Dixit et al ask in \cite{kzeta} whether a Ramanujan type identity for Hurwitz zeta exists and we answer this question positively. Ramanujan's identity gives an expression for the convolution of Riemann zeta at even arguments. In the same spirit it can be asked whether such an expression exists for the convolution of Riemann zeta at odd arguments, which was positively answered in \cite{dixit2021extended}:

	\begin{thm}\label{Atul_Herglotz}
		Let $\alpha$ and $\beta$ be two complex numbers such that $\Re(\alpha)>0, \Re(\beta)>0$ and $\alpha\beta=4\pi^2$. Let $\psi$ denote the digamma function. Then for $m \in \NN$, we have
		\begin{align*}&\left(-\beta\right)^{-m}\left\{2\gamma\zeta(2m+1)+\sum_{n=1}^{\infty}\dfrac{1}{n^{2m+1}}\left(\psi\left(\dfrac{\imath n\beta}{2\pi}\right)+\psi\left(-\dfrac{\imath n\beta}{2\pi}\right)\right)\right\}\\&+\alpha^{-m}\left\{2\gamma\zeta(2m+1)+\sum_{n=1}^{\infty}\dfrac{1}{n^{2m+1}}\left(\psi\left(\dfrac{\imath n\alpha}{2\pi}\right)+\psi\left(-\dfrac{\imath n\alpha}{2\pi}\right)\right)\right\}\\&=-2\sum_{k=1}^{m-1}\left(-1\right)^{k}\zeta(2k+1)\,\zeta(2m-2k+1)\alpha^{k-m}\beta^{-k}.\end{align*}
	\end{thm}
	The above identity can also be found in \cite{dirichlet}.
	
	It can be observed that 
	\begin{align*}  
		\frac{1}{2i\pi}\int_{(c)} \frac{\zeta(1-s) x^{-s}}{2\sin\left(\frac{\pi s}{2}\right)} ds &= \frac{\log(x)+\gamma}{\pi} -\frac{x^2}{\pi} \sum_{n=1}^{\infty} \frac{1}{n(x^2+n^2)}  \\
		&= \frac{1}{\pi} \left(\log(x)  - \frac{\psi(ix)+\psi(-ix)}{2}\right),
	\end{align*}
	where $c>1$. We generalize Theorem \ref{Atul_Herglotz} to its Hurwitz zeta analog involving another new kernel.
	
	\begin{defn} Let $x \in \RR^{+} , a \in \CC$ and $k \in \NN$. Define the \textit{odd Hurwitz kernel} by
		\begin{align}
			&\Phi(x,a;k) = \int_{(c)} \frac{\zeta(1-s,a)}{2\sin\left(\frac{\pi s}{2k}\right)} x^{-s} ds\\
			&= \frac{\log(x)+\gamma_{0}(a)}{\pi}-\frac{x^{2k}}{\pi}\sum_{n=0}^{\infty}\frac{1}{(n+a)((n+a)^{2k}+x^{2k})}
		\end{align}
	\end{defn}

	Finally we give a relation between kernels $\Phi(x,a;k)$ and $\Psi(x,a;b)$ in Theorem \ref{kernelrelation}. Throughout the remainder of the paper we assume $k \in \NN$ and $a,b \in \RR^{+}/\NN$ unless specified. For ease of notation we define
	\[\Psi\left(\frac{\alpha x}{\pi} ,a ;k\right)=\Psi_{\alpha}(x,a;k) \,\, \text{and} \,\, \Phi\left(\frac{\alpha x}{\pi} ,a ;k\right)=\Phi_{\alpha}(x,a;k). \]
	We prove our main results except Theorem \ref{thm1} and Theorem \ref{thm2} in two ways, first using complex analysis and then using tools borrowed from \cite{dirichlet}. For the convenience of the reader we include Dirichlet series setup and notations which can be found in \cite{dirichlet}.
	
	\subsection{Notations} For a sequence of non-zero complex numbers $\left\{ x_{n}\right\} $ and a sequence of associated complex \textbf{weights} $\left\{ a_{n}\right\} $,
	a Dirichlet series is defined as 
	\begin{equation}\label{notation1}
		\zeta_{x,a}\left(N\right)=\sum_{n=1}^{\infty}\frac{a_{n}}{x_{n}^{N}}.
	\end{equation}
	
	Note that we assume that the Dirichlet series is convergent for $N\geqslant 1.$ If the series diverges at $N=1$, but has a finite abscissa of convergence, we can obtain analogous results. We call  the function 
	\begin{equation}
		\label{gf}    
		\psi_{x,a}\left(z\right)=\sum_{N=1}^{\infty}\zeta_{x,a}\left(N\right)z^{N}
	\end{equation}
	associated to $\zeta_{x,a}$ the \textbf{zeta generating function}. Using  geometric series, we can check that the generating function $\psi_{x,a}\left(z\right)$ can be  expressed in terms
	of the weights $\left\{ a_{n}\right\} $ and zeros $\left\{ x_{n}\right\} $
	as follows 
	
	\begin{equation}\label{notation2}
		\psi_{x,a}\left(z\right)=\sum_{n=1}^{\infty}a_{n}\,\dfrac{z}{x_{n}-z}.
	\end{equation}

	Finally, we introduce the following notation: the modified sequence
	of weights $\left\{ a.\psi_{y,b}\right\} _{n\geqslant 1}$ is defined by 
	\[
	\left(a.\psi_{y,b}\right)_{n}=a_{n}\psi_{y,b}\left(x_{n}\right)
	\]
	so that the corresponding Dirichlet series is
	\[
	\zeta_{x,a.\psi_{y,b}}\left(N\right)=\sum_{n=1}^{\infty}\frac{a_{n}\psi_{y,b}\left(x_{n}\right)}{x_{n}^{N}}.
	\]
	Similarly, we denote as $\left(a.\psi_{y,b}.\psi_{z,c}\right)$
	the sequence defined by 
	\[
	\left(a.\psi_{y,b}.\psi_{z,c}\right)_{n}=a_{n}\psi_{y,b}\left(x_{n}\right)\psi_{z,c}\left(x_{n}\right)
	\]
	with associated Dirichlet series
	\[
	\zeta_{x,a.\psi_{y,b}.\psi_{z,c}}\left(N\right)=\sum_{n=1}^{\infty}\frac{a_{n}\psi_{y,b}\left(x_{n}\right)\psi_{z,c}\left(x_{n}\right)}{x_{n}^{N}}.
	\]
	Finally the convolution of two Dirichlet series is defined as 
	\[
	\left(\zeta_{y,b}*\zeta_{x,a}\right)\left(N+1\right)=\sum_{k=1}^{N}\zeta_{y,b}\left(k\right)\zeta_{x,a}\left(N+1-k\right)
	\]
	and the $n-$fold convolution of $n$ Dirichlet series $\zeta_{x^{\left(1\right)},a^{\left(1\right)}},\dots,\zeta_{x^{\left(n\right)},a^{\left(n\right)}}$
	as
	\[
	\left(\stackrel[i=1]{n}{*}\zeta_{x^{\left(i\right)},a^{\left(i\right)}}\right)\left(N+1\right)=\sum\zeta_{x^{\left(1\right)},a^{\left(1\right)}}\left(k_{1}\right)\dots\zeta_{x^{\left(n\right)},a^{\left(n\right)}}\left(k_{n}\right)
	\]
	where the sum is over the set of indices 
	\[
	\left\{ \left(k_{1},k_2,\ldots,k_{n}\right):1\leqslant k_{i}\leqslant N,\sum_{i=1}^{n}k_{i}=N+1\right\}.
	\]
	The proof of Theorem \ref{main_thm} below can be found in \cite{dirichlet}
	\begin{thm}
		\label{main_thm} For a set of $n\geqslant 2$ Dirichlet series $\left\{ \zeta_{x^{\left(i\right)},a^{\left(i\right)}}\right\} {}_{1\leqslant i\leqslant n}$,
		we have, evaluated at argument $N+1$ removed for clarity,
		\[
		\stackrel[i=1]{n}{*}\zeta_{x^{\left(i\right)},a^{\left(i\right)}}
		=\sum_{i=1}^{n}\zeta_{x^{\left(i\right)},a^{\left(i\right)}.\prod_{1\leqslant k\ne i\leqslant n}\psi_{x^{\left(k\right)}}}
		\]
		The special case $n=2$ reads
		\begin{align}
			\zeta_{y,b}*\zeta_{x,a}
			=\zeta_{x,a.\psi_{y}}
			+\zeta_{y,b.\psi_{x}}
			\label{eq:2terms}
		\end{align}
	\end{thm}
	For the sake of clarity, let us rephrase this identity (\ref{eq:2terms}) in a more explicit way:
	\[
	\left(\zeta_{y,b}*\zeta_{x,a}\right)\left(N+1\right)=\sum_{n=1}^{\infty}\left\{\frac{a_{n}\psi_{y}\left(x_{n}\right)}{x_{n}^{N+1}}+\frac{b_{n}\psi_{x}\left(y_{n}\right)}{y_{n}^{N+1}}\right\},
	\]
	
	\section{Main results}\label{mainresults}
	
	The \textit{Bernoulli polynomials} $B_n(x)$ are defined through their generating function as
	\begin{equation}
		\frac {te^{xt}}{e^{t}-1}=\sum _{n=0}^{\infty }B_{n}(x) \frac{t^n}{n!},
	\end{equation} 
	and the Bernoulli numbers $B_n$ are defined by $B_n=B_n(0)$. The \textit{Hurwitz zeta function} for complex variable $s$ with $\Re(s)>1$ and $a \in \CC /\ZZ_{\leq 0}$ is defined by 
	\begin{equation}
		\zeta (s,a)=\sum _{n=0}^{\infty }{\frac {1}{(n+a)^{s}}}.
	\end{equation}
	The Hurwitz zeta function has the integral representation \cite{apostol1998introduction}
	\begin{equation}
		\zeta(s,a)=\frac{1}{\Gamma(s)}\int_{0}^{\infty} \frac{x^{s-1}e^{-ax}}{1-e^{-x}} \mathrm{d}s,
	\end{equation}
	valid for $\Re(s)>1$. Moreover, it has an analytic continuation represented by the following contour integral \cite{whittekar}:
	
	\begin{equation}\label{eqn}
		\zeta (s,a)={\frac {-\Gamma (1-s)}{2\pi i}}\int _{C}{\frac {(-z)^{s-1}e^{-az}}{1-e^{-z}}}dz,
	\end{equation}
	where $C$ is the Hankel contour counterclockwise around the positive real axis and the principal branch for exponentiation is used.   
	The Laurent series expansion of Hurwitz zeta zentered at $s=1$ is given by
	\[ \zeta (s,a)={\frac {1}{s-1}}+\sum _{n=0}^{\infty }{\frac {(-1)^{n}}{n!}}\gamma _{n}(a)(s-1)^{n}\]
	where $\gamma_{n}(a)$ are generalized Stieltjes constants.
	
	The integral (\ref{eqn}) defines $\zeta(s, a)$ for all $s\in\CC$, with a single pole at $s = 1$ and corresponding residue $1$. Observing the analogy between Ramanujan's kernel and $\Psi(x,a;k)$ we derive a Ramanujan type formula involving $\Psi(x,a;k)$  in the theorem below.

	\begin{thm}\label{hurwitzramanujan} Let $\alpha,\beta\in\RR^+$ such that $\alpha\beta=\pi^2$, let $a,b \in \CC/\ZZ_{\leq 0}$ and $k,N \in \NN$. Then, we have
		\begin{align}\label{genhurwitz}
			&	\beta^{k\left(N+1\right)-1}\sum_{n=0}^{\infty}\dfrac{\Psi_{\alpha}(n+b,a;k)}{\left(n+b\right)^{2k\left(N+1\right)-1}} - \left(-1\right)^N \alpha^{k\left(N+1\right)-1}\sum_{n=0}^{\infty}\dfrac{\Psi_{\beta}(n+a,b;k)}{\left(n+a\right)^{2k\left(N+1\right)-1}}
			\nonumber\\ & =\dfrac{\left(-1\right)^N\alpha^{k\left(N+1\right)-1}\zeta(2k\left(N+1\right)-1,a)}{2k\cos\left(\frac{\pi\left(k-1\right)}{2k}\right)}-\frac{\beta^{k\left(N+1\right)-1}\zeta(2k\left(N+1\right)-1,b)}{2k\cos\left(\frac{\pi\left(k-1\right)}{2k}\right)} 
			\nonumber\\ & 
			+\sum_{p=0}^{N+1} \left(-1\right)^{p+1}\zeta(2kp,a)\,\zeta(2k\left(N-p+1\right),b)\,\alpha^{kp-1}\beta^{k(N+1-p)-1}.
		\end{align}

	\end{thm}
	Substituting the power series expansion of $\Psi_{\alpha}(x,a;k)$ in equation \eqref{genhurwitz} we get
	\begin{align}\label{2.1substitute}
		&\beta^{k(N+1)-1}\sum_{n=0}^{\infty}\frac{1}{(n+b)^{2k(N+1)-1}}\sum_{n=0}^{\infty} \frac{\alpha^{k-1}(n+b)^{2k-1}}{\alpha^{k}(n+b)^{2k}+\beta^{k}(n+a)^{2k}}\\
		&-(-1)^N\alpha^{k(N+1)-1}\sum_{n=0}^{\infty}\frac{1}{(n+a)^{2k(N+1)-1}}\sum_{n=0}^{\infty} \frac{\beta^{k-1}(n+a)^{2k-1}}{\beta^{k}(n+a)^{2k}+\alpha^{k}(n+b)^{2k}}\nonumber\\
		&=\sum_{p=1}^{N} \left(-1\right)^{p+1}\zeta(2kp,a)\,\zeta(2k\left(N-p+1\right),b)\,\alpha^{kp-1}\beta^{k(N+1-p)-1}\nonumber
	\end{align}
	When $N= 2q+1 , q \in \NN$ and $\alpha=\beta=\pi$ we have 
	\begin{align}
		& \sum_{n=0}^{\infty} \frac{\psi_{\pi}(n+b,a;k)}{(n+b)^{4k(q+1)-1}}+\sum_{n=0}^{\infty} \frac{\psi_{\pi}(n+a,b;k)}{(n+a)^{4k(q+1)-1}} \\
		&= \frac{1}{\pi}\sum_{p=0}^{2q+2} (-1)^{p+1}\zeta(2kp,a)\zeta(2k(2q-p+2),b)\nonumber  \\
		&-\frac{1}{2k\cos\left(\frac{\pi\left(k-1\right)}{2k}\right)}\left(\zeta(4k(q+1)-1,a)+\zeta(4k(q+1)-1,b)\right)\nonumber.
	\end{align}
	Since
	\[\Psi_{\pi}(x,a;k)=\Psi(x,a;k)=\frac{2a-1}{2\pi x}-\frac{1}{2k\cos\left(\frac{\pi\left(k-1\right)}{2k}\right)} + \frac{1}{\pi }\sum_{n=0}^{\infty}\frac{x^{2k-1}}{x^{2k}+(n+a)^{2k}},\]
	we have 
	\begin{align}
		&\frac{2a-1(\zeta(4k(q+1),a)+\zeta(4k(q+1),b))}{2}+\sum_{n=0}^{\infty}\frac{1}{(n+b)^{4k(q+1)}}\sum_{j=0}^{\infty}\frac{1}{1+\left(\frac{j+a}{n+b}\right)^{2k}}\\
		&+\sum_{n=0}^{\infty}\frac{1}{(n+a)^{4k(q+1)}}\sum_{j=0}^{\infty}\frac{1}{1+\left(\frac{j+b}{n+a}\right)^{2k}}=\sum_{p=0}^{2q+2} (-1)^{p+1}\zeta(2kp,a)\zeta(2k(2q-p+2),b)\nonumber,
	\end{align}
	which can be written as 
	\begin{align}
		&\sum_{n=0}^{\infty}\frac{1}{(n+b)^{4k(q+1)}}\sum_{j=0}^{\infty}\frac{1}{1+\left(\frac{j+a}{n+b}\right)^{2k}}+\sum_{n=0}^{\infty}\frac{1}{(n+a)^{4k(q+1)}}\sum_{j=0}^{\infty}\frac{1}{1+\left(\frac{j+b}{n+a}\right)^{2k}}\\
		&=\sum_{p=1}^{2q+1} (-1)^{p+1}\zeta(2kp,a)\zeta(2k(2q-p+2),b)\nonumber.
	\end{align}
	With $a=b$ we have 
	\begin{align}
		& \sum_{n=0}^{\infty}\frac{1}{(n+a)^{4k(q+1)}}\sum_{j=0}^{\infty} \frac{1}{1+\left(\frac{j+a}{n+a}\right)^{2k}} =\frac{1}{2}\sum_{p=0}^{2q+1} (-1)^{p+1}\zeta(2kp,a)\zeta(2k(2q-p+2),a) 
	\end{align}
	
	Setting $a=b$ in Theorem \ref{hurwitzramanujan} we deduce:
	\begin{cor}\label{2.2} Let $\alpha\beta=\pi^2$, $a \in \CC/\ZZ_{\leq 0}$ and $k,N \in \NN$. We have the identity
		
		\begin{align}\label{hurwitzzetta}
			&	\beta^{k\left(N+1\right)-1}\sum_{n=0}^{\infty}\dfrac{\Psi_{\alpha}(n+a,a;k)}{\left(n+a\right)^{2k\left(N+1\right)-1}} - \left(-1\right)^N \alpha^{k\left(N+1\right)-1}\sum_{n=0}^{\infty}\dfrac{\Psi_{\beta}(n+a,a;k)}{\left(n+a\right)^{2k\left(N+1\right)-1}}
			\nonumber\\ & =\dfrac{\left(-1\right)^N\alpha^{k\left(N+1\right)-1}\zeta(2k \left(N+1\right)-1,a)}{2k\cos\left(\frac{\pi\left(k-1\right)}{2k}\right)}-\frac{\beta^{k\left(N+1\right)-1}\zeta(2k\left(N+1\right)-1,a)}{2k\cos\left(\frac{\pi\left(k-1\right)}{2k}\right)} 
			\nonumber\\ & 
			+\sum_{p=0}^{N+1} \left(-1\right)^{p}\zeta(2kp,a)\,\zeta(2k\left(N-p+1\right),a)\,\alpha^{kp-1}\beta^{k(N+1-p)-1}.
		\end{align}
	\end{cor}

	As a corollary, when $a=1$ we have the following convolution of Riemann zeta functions:
	\begin{cor}
		Let $\alpha \beta = \pi^2$ 
		\begin{align}
			\label{riemannzeta}	
			&\beta^{k\left(N+1\right)-1}\sum_{n=1}^{\infty}\dfrac{\Psi_{\alpha}(n,1;k)}{\left(n\right)^{2k\left(N+1\right)-1}} - \left(-1\right)^N \alpha^{k\left(N+1\right)-1}\sum_{n=1}^{\infty}\dfrac{\Psi_{\beta}(n,1;k)}{\left(n\right)^{2k\left(N+1\right)-1}}
			\nonumber\\ & =\dfrac{\left(-1\right)^N\alpha^{k\left(N+1\right)-1}\zeta(2k\left(N+1\right)-1)}{2k\cos\left(\frac{\pi\left(k-1\right)}{2k}\right)}-\frac{\beta^{k\left(N+1\right)-1}\zeta(2k\left(N+1\right)-1)}{2k\cos\left(\frac{\pi\left(k-1\right)}{2k}\right)} 
			\nonumber\\ 
			& 
			+\sum_{p=0}^{N+1} \left(-1\right)^{p}\zeta(2kp)\,\zeta(2k\left(N-p+1\right))\,\alpha^{kp-1}\beta^{k(N+1-p)-1}.
		\end{align}
	\end{cor}
	When $k=1$ in Theorem \ref{hurwitzramanujan} we recover the identity
	\begin{align}
		&\sum_{k=1}^{N}\left(-\alpha\right)^{N+1-k}\beta^{k}\zeta\left(2k,a\right)\zeta\left(2N+2-2k,b\right)\nonumber \\
		&=\dfrac{\beta^{N+1}}{2}\,i\sqrt{\frac{\alpha}{\beta}}\sum_{n=0}^{\infty}\frac{1}{\left(a+n\right)^{2N+1}}\left[\psi\left(b+i\sqrt{\frac{\alpha}{\beta}}\left(a+n\right)\right)-\psi\left(b-i\sqrt{\frac{\alpha}{\beta}}\left(a+n\right)\right)\right]\\
		&+\frac{\left(-\alpha\right)^{N+1}}{2}\,i\sqrt{\frac{\beta}{\alpha}}\sum_{n=0}^{\infty}\frac{1}{\left(b+n\right)^{2N+1}}\left[\psi\left(a+i\sqrt{\frac{\beta}{\alpha}}\left(b+n\right)\right)-\psi\left(a-i\sqrt{\frac{\beta}{\alpha}}\left(b+n\right)\right)\right]\nonumber\label{eq:Hurwitz_main}\end{align}
	which has been recently proved in \cite{dirichlet} using Theorem \ref{main_thm}. Their method allows for straightforward generalizations of these transformations. Conversely, our method provides rigorous convergence guarantees, which their method does not provide. 
	
	Ramanujan's identity \ref{eq:Ramanujan_main} has the following special case, namely, for $\alpha,\beta>0$ with $\alpha\beta=\pi^2$,
	\begin{equation}\label{2.5}
		\alpha^m\sum_{n=1}^{\infty}\frac{n^{2m-1}}{e^{2\alpha n}-1}-(-\beta)^m\sum_{n=1}^{\infty}\frac{n^{2m-1}}{e^{2\beta n}-1}= (\alpha^m - (-\beta)^m)\frac{B_{2m}}{4m} \qquad (m > 1),
	\end{equation}
	
	which can be found in \cite{Berndt,Ramanujan}. The following theorem provides a generalization of equation \eqref{2.5} involving the Hurwitz kernel $\Psi(x,a;k)$.

	\begin{thm}\label{thm1}
		Let $\alpha,\beta\in\mathbb{R}^+$ such that $\alpha\beta=\pi^2$ and $B_{j}(x)$ denote Bernoulli polynomials. Then, the following identity holds
		\begin{align}
			&\alpha^{km+1}\sum_{n=0}^{\infty}\left(n+b\right)^{2km+1}\left[\Psi_{\alpha}(n+b,a;k) - \sum_{p=1}^{m} \frac{B_{2kp+1}(a)}{2kp+1}\left(\frac{\pi}{\alpha}\right)^{2kp+1}\right]
			\\&+\left(-1\right)^m\beta^{km+1} \sum_{n=0}^{\infty}\left(n+a\right)^{2km+1}\left[\Psi_{\beta}(n+a,b;k) - \sum_{p=1}^{m} \frac{B_{2kp+1}(b)}{2kp+1}\left(\frac{\pi}{\beta}\right)^{2kp+1}\right] 
			\nonumber\\&= \frac{\left(-1\right)^{m+1}\beta^{km+1}B_{2km+2}(a)}{4k(km+1)\cos\left(\frac{\pi\left(k-1\right)}{2k}\right)}-\frac{\alpha^{km+1}B_{2km+2}(b)}{4k(km+1)\cos\left(\frac{\pi\left(k-1\right)}{2k}\right)}\nonumber\\
			&	+\sum_{p=0}^{m}(-1)^p \frac{B_{2kp+1}(a)\beta^{kp}B_{2k\left(m-p\right)+1}(b)\alpha^{k\left(m-p\right)}}{\left(2kp+1\right)\left(2k\left(m-p\right)+1\right)}\nonumber.
		\end{align}
	\end{thm}
	
	As a special case when $a=b=1$ we have:
	
	\begin{prop}\label{prop}
		Let $\alpha,\beta>0$ and $\alpha\beta=\pi^2$. The following identity holds 
		\begin{align}
			&(\alpha)^{km+1} \sum_{n=1}^{\infty} n^{2km+1}\Psi_{\alpha}(n,1;k)+(-1)^{m}\beta^{km+1}\sum_{n=1}^{\infty} n^{2km+1}\Psi_{\beta}(n,1;k) \\ 
			&=\frac{1}{4k(km+1)\cos\left(\frac{\pi(k-1)}{2k}\right)}\left(\alpha^{km+1}B_{2km+2}-(-1)^{m+1}\beta^{km+1}B_{2km+2}\right)\nonumber.
		\end{align}
	\end{prop}
	
	Moreover, plugging in $\alpha=\beta=\pi$ with $m=2p$ in Proposition \ref{prop} produces
	\begin{equation}\label{R}
		\sum_{n=1}^{\infty} \dfrac{\Psi_{\pi}(n,1;k)}{n^{-4kp-1} } =\frac{B_{4kp+2}}{4k\left(2kp+1\right)\cos\left(\frac{\pi\left(k-1\right)}{2k}\right)},
	\end{equation}
	which is a generalization of Glaisher's famous identity \cite{Glashier}
	\begin{equation}
		\sum_{n=1}^{\infty}\frac{n^{4m+1}}{e^{2\pi n}-1} =  \frac{B_{4m+2}}{2(4m+2)}.
	\end{equation}
	When we put $k=1$ in Proposition \ref{prop} we recover equation \eqref{2.5}.
	
	As a companion of Theorem $\ref{thm1}$ we have the following identity for $\Phi(x,a;k)$
	
	\begin{thm}\label{thm2}
		Let $\alpha,\beta>0$ and $\alpha\beta=\pi^2$. The following identity holds
		\begin{align}
			&\alpha^{km}\sum_{n=0}^{\infty}\left(n+b\right)^{2km-1}\left[\Phi_{\alpha}(n+b,a;k)-\sum_{p=1}^{m} \frac{B_{2kp}(a)}{2kp}\left(\frac{\pi}{\alpha}\right)^{2kp}\right]\\ &+(-1)^m\beta^{km} \sum_{n=0}^{\infty}\left(n+a\right)^{2km-1}\left[\Phi_{\beta}(n+a,b;k)-\sum_{p=1}^{m} \frac{B_{2kp}(b)}{2kp}\left(\frac{\pi}{\beta}\right)^{2kp}\right] \nonumber\\
			&=\frac{-\alpha^{km}}{\pi}\left(\frac{B_{2km}(b)}{2km}\left(\log\left(\frac{\alpha}{\pi}\right) + \gamma_{0}(a)\right) +\frac{\partial}{\partial s}\zeta(s-2km+1,b) |_{s=0} \right)\nonumber\\
			&+\frac{(-1)^{m+1}\beta^{km}}{\pi}\left(\frac{B_{2km}(a)}{2km}\left(\log\left(\frac{\beta}{\pi}\right) +\gamma_0(b)\right)+\frac{\partial}{\partial s} \zeta(s-2km+1,a)\big |_{s=2km}\right)\nonumber\\
			& +\sum_{p=1}^{m-1}(-1)^p\ \frac{\beta^{kp}B_{2kp}(a)\alpha^{k(m-p)}B_{2k\left(m-p\right)}(b)}{\pi\left(2kp\right)\left(2k\left(m-p\right)\right)}	\nonumber		
		\end{align}
	\end{thm}

	We now give Hurwitz zeta generalization of Theorem \ref{Atul_Herglotz} involving the kernel $\Phi(x,a;k)$.
	\begin{thm}\label{Atul_Herglotz_1}
		Let $\alpha,\beta>0$ and $\alpha\beta=\pi^2$. The following identity holds 
		\begin{align}
			&\beta^{km}\sum_{n=0}^{\infty}\frac{\Phi_{\alpha}(n+b,a;k)}{(n+b)^{2km+1}}+(-\alpha^k)^m\sum_{n=0}^{\infty}\frac{\Phi_{\beta}(n+a,b;k)}{(n+a)^{2m+1}}\\
			&=\sum_{i=1}^{m-1}(-1)^i\frac{\alpha^{ki}\zeta(2ki+1,a)\beta^{k(m-i)}\zeta(2k(m-i)+1,b)}{\pi}\nonumber\\
			&+\frac{\beta^{km}}{\pi}\left(\zeta(2km+1,b)\left(\log\left(\frac{\alpha}{\pi}\right)+\gamma_{0}(a)\right)-\frac{\partial}{\partial s}\zeta(2km+1+s,b) \big |_{s=0}\right)\nonumber\\
			&+\frac{(-1)^m\alpha^{km}}{\pi} \left(\zeta(2km+1,a)\left(\gamma_{0}(b)-\log\left(\frac{\alpha}{\pi}\right)\right)-\frac{\partial}{\partial s}\zeta(2km+1+s,a) \big |_{s=0} \right) \nonumber
		\end{align}
	\end{thm}
	We recover Theorem \ref{Atul_Herglotz} when we put $a=b=k=1$ in Theorem \ref{Atul_Herglotz_1}.
	
	Substituting the power series expansion of $\Phi_{\alpha}(n+b,a;k)$ and $\Phi_{\beta}(n+a,b;k)$ we have
	\begin{align*}
		&\sum_{n=0}^{\infty} \frac{\Phi_{\alpha}(n+b,a;k)}{(n+b)^{2km+1}} \\
		& =\sum_{n=0}^{\infty} \frac{1}{(n+b)^{2km+1}}\left(\frac{\log\left(\frac{\alpha (n+b)}{\pi}\right) +\gamma_{0}(a)}{\pi}-\frac{\alpha^k(n+b)^{2k}}{\pi}\sum_{i=0}^{\infty}\frac{1}{(i+a)(\beta^k(i+a)^{2k} +\alpha^k(n+b)^{2k})}\right)\\
		&\sum_{n=0}^{\infty} \frac{\Phi_{\beta}(n+a,b;k)}{(n+a)^{2km+1}} \\
		& =\sum_{n=0}^{\infty} \frac{1}{(n+a)^{2km+1}}\left(\frac{\log\left(\frac{\beta (n+a)}{\pi}\right) +\gamma_{0}(b)}{\pi}-\frac{\beta^k(n+a)^{2k}}{\pi}\sum_{i=0}^{\infty}\frac{1}{(i+b)(\beta^k(i+b)^{2k} +\alpha^k(n+a)^{2k})}\right)\\
	\end{align*}
	Thus the identity can be written as
	\begin{align}\label{eqtransform_2}
		&-\pi^{2k-1}\beta^{k(m-1)}\sum_{n=0}^{\infty}\frac{1}{(n+b)^{2k(m-1)+1}}\sum_{i=0}^{\infty}\frac{1}{(i+a)(\beta^k(i+a)^{2k} +\alpha^k(n+b)^{2k})}\\
		& -\pi^{2k-1}(-\alpha)^{k(m-1)}\sum_{n=0}^{\infty}\frac{1}{(n+a)^{2k(m-1)+1}}\sum_{i=0}^{\infty}\frac{1}{(i+b)(\beta^k(i+b)^{2k} +\alpha^k(n+a)^{2k})}\nonumber\\
		&= \sum_{i=1}^{m-1}(-1)^i\frac{\alpha^{ki}\zeta(2ki+1,a)\beta^{k(m-i)}\zeta(2k(m-i)+1,b)}{\pi}\nonumber
	\end{align}

	In \cite{dixit2020generalized} , Dixit et al gave an expression for convolution of riemann zeta function at odd and even values, namely, for $\alpha,\beta>0$ with $\alpha \beta = \pi^2$ and $m \in \NN$
	
	\begin{equation}\label{atuloddeven}
		\beta^{-\left(m-\frac{1}{2}\right)}\left\{\frac{1}{2}\zeta(2m) +\sum_{n=0}^{\infty} \frac{n^{-2m}}{e^{2n \beta }-1}\right\}-\sum_{k=0}^{m-1}(-1)^{k+1}\frac{\zeta(2k)\zeta(2m-2k+1)}{\pi^{2k}}\beta^{2k-m-\frac{1}{2}} \end{equation}
	\begin{equation}= (-1)^{m+1}\alpha^{-\left(m-\frac{1}{2}\right)}\left\{\frac{\gamma}{\pi}\zeta(2m)+\frac{1}{2\pi} \sum_{n=1}^{\infty}n^{-2m}\left(\psi\left(\frac{in\alpha}{\pi}\right)+\psi\left(\frac{-in\alpha}{\pi}\right)\right)\right\}.\nonumber\end{equation}

	Equation \eqref{atuloddeven} is  natural companion to Ramanujan's identity and can be interpreted as a relation between $\Phi_{\alpha}(n,1;1)$ and $\Psi_{\beta}(n,1;1)$. The following t``heorem provides extension of this relation to Hurwitz zeta function and containing the kernes $\Phi_{\alpha}(x,a;k)$ and $\Psi_{\beta}(x,a;k)$ 
	
	\begin{thm}\label{kernelrelation}
		Let $\alpha,\beta>0$ and $\alpha \beta = \pi^2$. The following identity holds:
		\begin{align}
			&\beta^{km}\sum_{n=0}^{\infty} \frac{\Phi_{\alpha}(n+b,a;k)}{(n+b)^{2km}}+\pi(-1)^{m-1}\alpha^{km-1}\sum_{n=0}^{\infty}\frac{\Psi_{\beta}(n+a,b;k)}{(n+a)^{2km}}\\
			&=\frac{\beta^{km}}{\pi}\left(\zeta(2km,b)\left(\log\left(\frac{\alpha}{\pi}\right) +\gamma_0(a)\right)-\frac{\partial}{\partial s}\zeta(2km+s,b)\big |_{s=0}\right)+(-1)^m\pi\alpha^{km-1}\frac{\zeta(2km,a)}{2k \sin\left(\frac{\pi}{2k}\right)}\nonumber\\
			&+\sum_{i=1}^{m}\frac{(-1)^i\alpha^{ki}\zeta(2ki+1,a)\beta^{k(m-i)}\zeta(2k(m-i),b)}{\pi}\nonumber
		\end{align}
	\end{thm}
	When we put $a=b=k=1$ in Theorem \ref{kernelrelation} we recover equation \eqref{atuloddeven}. 
	Substituting the power series expansion of $\Phi_{\alpha}(n+b,a;k) ,\Psi_{\beta}(n+a,b;k) $ we have 
	\begin{align*}
		&\sum_{n=0}^{\infty} \frac{\Phi_{\alpha}(n+b,a;k)}{(n+b)^{2km}} \\
		& =\sum_{n=0}^{\infty} \frac{1}{(n+b)^{2km}}\left(\frac{\log\left(\frac{\alpha (n+b)}{\pi}\right) +\gamma_{0}(a)}{\pi}-\frac{\alpha^k(n+b)^{2k}}{\pi}\sum_{i=0}^{\infty}\frac{1}{(i+a)(\beta^k(i+a)^{2k} +\alpha^k(n+b)^{2k})}\right)\\
		& \sum_{n=0}^{\infty}\frac{\Psi_{\beta}(n+a,b;k)}{(n+a)^{2km}} \\
		& = \frac{2b-1}{2\beta}\zeta(2km+1,a)-\frac{\zeta(2km,a)}{2k\cos
			\left(\frac{\pi}{2k}(k-1)\right)}+\sum_{n=0}^{\infty}\frac{1}{(n+a)^{2k(m-1)+1}}\sum_{i=0}^{\infty}\frac{\beta^{k-1}}{\alpha^k(n+a)^{2k}+\beta^k(i+b)^{2k}}
	\end{align*}
	Thus the identity can be written as 
	\begin{align}\label{eqtransfrom_1}
		& \pi^{2k-1}(-1)^{m-1}\alpha^{k(m-1)}\sum_{n=0}^{\infty}\frac{1}{(n+a)^{2k(m-1)+1}}\sum_{i=0}^{\infty}\frac{1}{\beta^k(n+a)^{2k}+\alpha^k(i+b)^{2k}}\\
		&-\pi^{2k-1}\beta^{k(m-1)}\sum_{n=0}^{\infty} \frac{1}{(n+b)^{2k(m-1)}}\sum_{i=0}^{\infty}\frac{1}{(i+a)(\beta^k(i+a)^{2k} +\alpha^k(n+b)^{2k})}\nonumber\\
		&+\sum_{i=1}^{m-1}\frac{(-1)^i\alpha^{ki}\zeta(2ki+1,a)\beta^{k(m-i)}\zeta(2k(m-i),b)}{\pi}\nonumber\end{align}

	\section{Proofs}

	Before proving these results we state an important inequality which we will use throughout the proofs. Stirling’s formula
	for the gamma function on a vertical strip states that for $a  \leq \sigma  \leq b$ and $|t| \geq 1$,
	\begin{equation}\label{gamma}|\Gamma(\sigma+it)| = (2\pi)^{\frac{1}{2}}|t|^{\sigma-\frac{1}{2}}e^{-\frac{\pi}{2}|t|}\left(1+O\left(\frac{1}{|t|}\right)\right).\end{equation}
	The reflection formula for the gamma function is 
	$\Gamma (1-z)\Gamma (z)={\frac {\pi }{\sin \pi z}}, z\not \in \mathbb {Z}$
	Thus as $\Im(s)\to \infty$ using the reflection formula and \eqref{gamma} we have the inequality
	\begin{equation}\label{sine}\frac{1}{\left|\sin\left(\frac{\pi s}{2k}\right)\right|} = 2\exp\left(\frac{-\pi}{2}\left|\frac{\Im(s)}{k}\right|\right)\left(1+O\left(\frac{1}{|\Im(s)|}\right)\right). \end{equation}
	
	A variant of the reflection formula for the gamma function is $\Gamma\left(\frac{1}{2}+z\right)\Gamma\left(\frac{1}{2}-z\right)=\frac{\pi}{\cos(\pi z)} , z\notin \ZZ-\frac{1}{2}$. Thus as $\Im(s)\to \infty$ using the variant of reflection formula and \eqref{gamma} we have the inequality
	\begin{equation}\label{cosine}\frac{1}{\left|\cos\left(\frac{\pi(s+k-1)}{2k}\right)\right|}= 2\exp\left(\frac{-\pi}{2}\left|\frac{\Im(s)}{k}\right|\right)\left(1+O\left(\frac{1}{|\Im(s)|}\right)\right)\end{equation}

	\subsection{First Proof of Theorem \ref{hurwitzramanujan}}
	
	Let $\Psi_{\alpha}(x) = \Psi\left(\frac{\alpha x}{\pi}\right)$. On account of absolute convergence of the Hurwitz zeta function for $\Re(s)>1$ we have 
	\[\sum_{n=0}^{\infty}(n+b)^{-2kN-2k+1}\Psi_{\alpha}\left(n+b,a\right) = \int_{(c)} \frac{\zeta(1-s,a)\zeta(2kN+2k-1+s,b)}{2k\cos\left(\frac{\pi(s+k-1)}{2k}\right)} \left(\frac{\alpha}{\pi}\right)^{-s} ds. \]
	We now evaluate this integral by shifting the line of integration. Consider rectangular the contour determined by the line segments $[c-iT,c+iT],[c+iT,d+iT],[d+iT,d-iT],[d-iT,c-iT]$ where $d = -c-2kN-2k+2$. Inside this domain, the integrand has simple poles at $0,-2kN-2k+2$ due to $\zeta(1-s,a)$ and $\zeta(2kN+2k-1+s,a)$ respectively. It also has simple poles at the integers $-2kp+1$ where $p \in \{0,1,\ldots,N+1\}$ due to the cosine term in the denominator. The residues at these poles are 
	\begin{align*}
		&R_0 = \frac{-\zeta(2kN+2k-1,b)}{2k\cos\left(\frac{\pi(k-1)}{2k}\right)}, \\
		&R_{-2kN-2k+2} = \left(\frac{\alpha}{\pi}\right)^{2kN+2k-2}\frac{\zeta(2kN+2k-1,a)}{2k\cos\left(\frac{\pi(-2kN-k+1)}{2k}\right)},\\
		&R_{-2kp+1} = (-1)^{p}\left(\frac{\alpha}{\pi}\right)^{2kp-1}\zeta(2kp,a)\zeta(2k(N+1-p),b).
	\end{align*}
	Thus by Cauchy's residue formula we have 
	\begin{align*}
		&\frac{1}{2i\pi}\left[\int_{c-iT}^{c+iT}+\int_{c+iT}^{d+iT}+\int_{d+iT}^{d-iT}+\int_{d-iT}^{c-iT}\right] \frac{\zeta(1-s,a)\zeta(2kN+2k-1+s,b)}{2k\cos\left(\frac{\pi(s+k-1)}{2k}\right)} \left(\frac{\alpha}{\pi}\right)^{-s} \mathrm{d}s\\
		&= R_{0}+R_{-2kN-2k+2}+\sum_{p=0}^{N+1} R_{-2kp+1}.
	\end{align*}
	From elementary bounds on the Hurwitz zeta function and Equation \eqref{cosine}, it can be seen that as $T \to \infty$, the integrals along horizontal segments tend to zero. Under the change of variables $ s \to -s-2kN-2k+2$ we have 
	\[\int_{(d)} \frac{\zeta(1-s,a)\zeta(2kN+2k-1+s,b)}{2k\cos\left(\frac{\pi(s+k-1)}{2k}\right)} \left(\frac{\alpha}{\pi}\right)^{-s} \]\[= (-1)^N\left(\frac{\alpha}{\pi}\right)^{2kN+2k-2}\int_{(c)}\frac{\zeta(1-s,b)\zeta(2kN+2k-1+s,a)}{2k\cos\left(\frac{\pi(s+k-1)}{2k}\right)} \left(\frac{\beta}{\pi}\right)^{-s}, \] 
	
	which proves Theorem \ref{hurwitzramanujan}.
	\subsection{Second Proof of Theorem \ref{hurwitzramanujan}}
	Choose
	\[a_n=b_n=1 \,\, , x_{n}=\frac{-(n+a)^{2k}}{\alpha^k} \,\,\text{and} \,\, y_n = \frac{(n+b)^{2k}}{\beta^{k}}\]
	which has the corrosponding zeta function and zeta generating function as 
	\[\zeta_{x,a}(q) = (-1)^q\alpha^{kq}\zeta(2kq,a) \,\, , \,\, \psi_{x,a}(z) = \sum_{n=0}^{\infty}\frac{\alpha^kz}{(n+a)^{2k}-\beta^kz}\]
	Using Theorem \ref{main_thm} we have 
	\begin{align}
		&\sum_{p=1}^{N}(-1)^p\zeta(2kp,a)\alpha^{kp}\zeta(2k(N+1-p),b)\beta^{k(N+1-p)}\\
		&=(-1)^N\alpha^{k(N+1)-1}\sum_{n=0}^{\infty}\frac{1}{(n+a)^{2k(N+1)-1}}\sum_{n=0}^{\infty} \frac{\beta^{k-1}(n+a)^{2k-1}}{\beta^{k}(n+a)^{2k}+\alpha^{k}(n+b)^{2k}}\nonumber\\
		&-\beta^{k(N+1)-1}\sum_{n=0}^{\infty}\frac{1}{(n+b)^{2k(N+1)-1}}\sum_{n=0}^{\infty} \frac{\alpha^{k-1}(n+b)^{2k-1}}{\alpha^{k}(n+b)^{2k}+\beta^{k}(n+a)^{2k}}\nonumber	
	\end{align}
	which can be seen equivalent to equation \eqref{2.1substitute}.

	\subsection{Proof of Theorem \ref{thm1}}
	On account of the absolute convergence of $\zeta(s,a)$ for $\Re(s)>1$ we have 
	\begin{align*}
		&\sum_{n=0}^{\infty}\left(n+b\right)^{2km+1}\left[\Psi_{\alpha}(n+b,a;k) - \sum_{p=1}^{m} \frac{B_{2kp+1}(a)}{2kp+1}\left(\frac{\pi}{\alpha}\right)^{2kp+1}\right] \\
		&=\frac{1}{2i\pi}\int_{(d)} \frac{\zeta(1-s,a)\zeta(s-2km-1,b)}{2k\cos\left(\frac{\pi(s+k-1)}{2k}\right)}\left(\frac{\alpha}{\pi}\right)^{-s}\mathrm{d}s,
	\end{align*}
	where $2km+2<d<2km+3$ since 
	$$\int_{(c)}\frac{\zeta(1-s,a)}{2k\cos\left(\frac{\pi(s+k-1)}{2k}\right)}x^{-s}\mathrm{d}s = \int_{(d)}\frac{\zeta(1-s,a)}{2k\cos\left(\frac{\pi(s+k-1)}{2k}\right)}x^{-s}\mathrm{d}s- \sum_{p=1}^{m} \frac{B_{2kp+1}(a)}{2kp+1}x^{-2kp-1}.$$
	
	We now evaluate this integral by shifting the line of integration and using the Cauchy's Residue Theorem. Consider the rectangular contour determined by the line segments $[d-iT,d+iT],[d+iT,e+iT],[e+iT,e-iT],[e-iT,d-iT]$ where $e = 2km+2-d$. Inside this domain, the integrand has simple poles at $0,2km+2$ due to $\zeta(1-s,a)$ and $\zeta(s-2km-2,b)$ respectively. It also has simple poles at the integers $2kp+1$ where $p \in \{0,1,\ldots,m\}$ due to cosine term in the denominator. The residues at these poles are 
	\begin{align*}
		& R_{0} = -\frac{B_{2km+2}(b)}{2k(2km+2)\cos\left(\frac{\pi\left(k-1\right)}{2k}\right)}, \\
		& R_{2km+2} = (-1)^{m+1}\left(\frac{\pi}{\alpha}\right)^{2km+2}\frac{\alpha^{km+1}B_{2km+2}(b)}{2k(2km+2)\cos\left(\frac{\pi\left(k-1\right)}{2k}\right)},\\
		& R_{2kp+1} = (-1)^p\left(\frac{\pi}{\alpha}\right)^{2kp+1} \frac{B_{2kp+1}(a)B_{2k\left(m-p\right)+1}(b)}{\pi\left(2kp+1\right)\left(2k\left(m-p\right)+1\right)}.
	\end{align*}
	Thus by Cauchy's residue formula we have 
	\begin{align*}
		&\frac{1}{2i\pi}\left[\int_{d-iT}^{d+iT}+\int_{d+iT}^{e+iT}+\int_{e+iT}^{e-iT}+\int_{e-iT}^{d-iT}\right] \frac{\zeta(1-s,a)\zeta(s-2km-1,b)}{2k\cos\left(\frac{\pi(s+k-1)}{2k}\right)}\left(\frac{\alpha}{\pi}\right)^{-s}\mathrm{d}s\\
		&= R_{0}+R_{2km+2}+\sum_{p=0}^{m} R_{2kp+1}.
	\end{align*}
	From elementary bounds on the Hurwitz zeta and Equation \ref{cosine}, it can be seen that as $T \to \infty$, the integrals along horizontal segments tend to zero. Under the change of variables $ s \to 2km+2-s$ we have 
	\begin{align*}
		&\int_{(e)}\frac{\zeta(1-s,a)\zeta(s-2km-1,b)}{2k\cos\left(\frac{\pi(s+k-1)}{2k}\right)}\left(\frac{\alpha}{\pi}\right)^{-s}\mathrm{d}s \\
		& = \left(\frac{\alpha}{\pi}\right)^{-2km-2}(-1)^m\int_{(d)}\frac{\zeta(1-s,b)\zeta(s-2km-1,a)}{2k\cos\left(\frac{\pi(s+k-1)}{2k}\right)}\left(\frac{\beta}{\pi}\right)^{-s}\mathrm{d}s
	\end{align*}
	which proves Theorem \ref{thm1}.
	
	\subsection{Proof of Theorem \ref{thm2}} 
	On account of the absolute convergence of $\zeta(s,a)$ for $\Re(s)>1$ we have 
	\begin{align*}
		&\sum_{n=0}^{\infty}\left(n+b\right)^{2km-1}\left[\Phi_{\alpha}(n+b,a;k) - \sum_{p=1}^{m} \frac{B_{2kp}(a)}{2kp}\left(\frac{\pi}{\alpha}\right)^{2kp}\right] \\
		&=\frac{1}{2i\pi}\int_{(d)} \frac{\zeta(1-s,a)\zeta(s-2km+1,b)}{2k\sin\left(\frac{\pi s}{2k}\right)}\left(\frac{\alpha}{\pi}\right)^{-s}\mathrm{d}s,
	\end{align*}
	where $2km<d<2km+1$ since 
	
	$$\int_{(c)}\frac{\zeta(1-s,a)}{2k\sin\left(\frac{\pi s}{2k}\right)}x^{-s}\mathrm{d}s = \int_{(d)}\frac{\zeta(1-s,a)}{2k\sin\left(\frac{\pi s}{2k}\right)}x^{-s}\mathrm{d}s- \sum_{p=1}^{m} \frac{B_{2kp}(a)}{2kp}x^{-2kp}.$$

	We now evaluate this integral by shifting the line of integration and using the Residue Theorem. Consider the contour determined by the line segments $[d-iT,d+iT],[d+iT,e+iT],[e+iT,e-iT],[e-iT,d-iT]$ where $e = 2km-d$. The integrand has poles of order two at $0,2km$. It also has simple poles at the integers $2kp$ where $p \in \{1,\ldots,m-1\}$ due to cosine term in the denominator. The residues at these poles are 
	\begin{align*}
		& R_{0} = \frac{-1}{\pi}\left(\frac{B_{2km}(b)}{2km}\left(\log\left(\frac{\alpha}{\pi}\right) + \gamma_{0}(a)\right) +\frac{\partial}{\partial s}\zeta(s-2km+1,b) |_{s=0} \right) \\
		& R_{2km} =\frac{(-1)^m}{\pi}\left(\frac{\pi}{\alpha}\right)^{2km}\left(\frac{B_{2km}(a)}{2km}\left(\log\left(\frac{\alpha}{\pi}\right) -\gamma_0(b)\right)-\frac{\partial}{\partial s} \zeta(s-2km+1,a)\big |_{s=2km}\right) \\
		& R_{2kp} = (-1)^p\left(\frac{\pi}{\alpha}\right)^{2kp} \frac{B_{2kp}(a)B_{2k\left(m-p\right)}(b)}{\pi\left(2kp\right)\left(2k\left(m-p\right)\right)}.
	\end{align*}
	Thus we have 
	\begin{align*}
		&\frac{1}{2i\pi}\left[\int_{d-iT}^{d+iT}+\int_{d+iT}^{e+iT}+\int_{e+iT}^{e-iT}+\int_{e-iT}^{d-iT}\right] \frac{\zeta(1-s,a)\zeta(s-2km+1,b)}{2k\sin\left(\frac{\pi s}{2k}\right)}\left(\frac{\alpha}{\pi}\right)^{-s}\mathrm{d}s\\
		&= R_{0}+R_{2km}+\sum_{p=1}^{m-1} R_{2kp}.
	\end{align*}
	From elementary bounds on the Hurwitz zeta and Equation \ref{cosine}, it can be seen that as $T \to \infty$, the integrals along horizontal segments tend to zero. Under the change of variables $ s \to 2km-s$ we have 
	\begin{align*}
		&\int_{(e)}\frac{\zeta(1-s,a)\zeta(s-2km-1,b)}{2k\sin\left(\frac{\pi s}{2k}\right)}\left(\frac{\alpha}{\pi}\right)^{-s}\mathrm{d}s \\
		& = \left(\frac{\alpha}{\pi}\right)^{-2km}(-1)^{m-1}\int_{(d)}\frac{\zeta(1-s,b)\zeta(s-2km-1,a)}{2k\sin\left(\frac{\pi s}{2k}\right)}\left(\frac{\beta}{\pi}\right)^{-s}\mathrm{d}s
	\end{align*}
	which proves Theorem \ref{thm2}.
	
	\subsection{First Proof of Theorem \ref{Atul_Herglotz_1}}
	On account of absolute convergence of $\zeta(s,b)$ for $\Re(s)>1$ we have 
	\begin{equation}
		\sum_{n=0}^{\infty}\frac{\Phi_{\alpha}(n+b,a;k)}{(n+b)^{2km+1}}=\int_{(c)}\frac{\zeta(1-s,a)\zeta(2km+1+s,b)}{2\sin\left(\frac{\pi s}{2k}\right)} \left(\frac{\alpha}{\pi}\right)^{-s} ds
	\end{equation}
	where $1<c<2$.
	We now evaluate this integral by shifting the line of integration and using Residue Theorem. Consider the recangular contour determined by the line segments $[c-iT,c+iT],[c+iT,d+iT],[d+iT,d-iT],[d-iT,c-iT]$ where $d = -2km-c$. Inside this domain, the integrand has simple poles at $0,-2km$ due to $\zeta(1-s,a)$ and $\zeta(s-2m-1,b)$ respectively. It also has simple poles at the integers $-2ki$ where $i \in \{1,\ldots,m-1\}$ due to sine term in the denominator. The residues at these poles are 
	\begin{align*}
		& R_{0} = \frac{1}{\pi}\left(\zeta(2km+1,b)\left(\log\left(\frac{\alpha}{\pi}\right)+\gamma_{0}(a)\right)-\frac{\partial}{\partial s}\zeta(2km+1+s,b) \big |_{s=0}\right) \\
		& R_{-2km}=\frac{(-1)^m}{\pi}\left(\frac{\alpha}{\pi}\right)^{2km} \left(\zeta(2km+1,a)\left(\gamma_{0}(b)-\log\left(\frac{\alpha}{\pi}\right)\right)+\frac{\partial}{\partial s}\zeta(2km+1+s,a) \big |_{s=0} \right) \\
		& R_{-2i}=(-1)^i \left(\frac{\alpha}{\pi}\right)^{2ki}\frac{\zeta(2ki+1,a)\zeta(2k(m-i)+1,b)}{\pi}
	\end{align*}
	Thus we have 
	\begin{align*}
		&\frac{1}{2i\pi}\left[\int_{c-iT}^{c+iT}+\int_{c+iT}^{d+iT}+\int_{d+iT}^{d-iT}+\int_{d-iT}^{c-iT}\right] \frac{\zeta(1-s,a)\zeta(s+2km+1,b)}{2\sin\left(\frac{\pi s}{2k}\right)}\left(\frac{\alpha}{\pi}\right)^{-s}\mathrm{d}s\\
		&= R_{0}+R_{-2km}+\sum_{i=1}^{m-1} R_{-2ki}
	\end{align*}
	From elementary bounds on Hurwitz zeta and \eqref{cosine}, it can be seen that as $T \to \infty$, the integrals along horizontal segments tend to zero. Under the change of variables $ s \to -2km-s$ we have 
	\begin{align*}
		&\int_{(d)} \frac{\zeta(1-s,a)\zeta(s+2km+1,b)}{2\sin\left(\frac{\pi s}{2}\right)}\left(\frac{\alpha}{\pi}\right)^{-s}\mathrm{d}s \\
		& = \left(\frac{\alpha}{\pi}\right)^{2km}(-1)^m\int_{(c)} \frac{\zeta(1-s,b)\zeta(s-2km-1,a)}{2\sin\left(\frac{\pi s}{2k}\right)}\left(\frac{\beta}{\pi}\right)^{-s}\mathrm{d}s
	\end{align*}
	which proves Theorem \ref{Atul_Herglotz_1}.
	\subsection{Second Proof of Theorem \ref{Atul_Herglotz_1}}
	Choose
	\[a_n=\frac{1}{n+a} \,\, x_n = -\frac{(n+a)^{2k}}{\alpha^k} \,\, b_n = \frac{1}{n+b} \,\, \text{and} \,\, y_n = \frac{(n+b)^{2k}}{\beta^k}\]
	so that the corrosponding zeta functions are \[\zeta_{x,a}(i) = (-1)^i\alpha^{ki}\zeta(2ki+1,a) \,\, , \,\, \zeta_{y,b}(m-i) = \beta^{k(m-i)}\zeta(2k(m-i)+1,b)\]
	and corrosponding zeta generating functions are 
	\[\psi_{x,a}(z) = -\sum_{n=0}^{\infty}\frac{\alpha^k z}{(n+a)((n+a)^{2k}+\alpha^k z)} \,\, , \,\, \psi_{y,b}(z)=\sum_{n=0}^{\infty}\frac{\beta^k z}{(n+b)((n+b)^{2k}-\beta^k z)}\]
	Using Theorem \ref{main_thm} we have
	\begin{align}
		&= \sum_{i=1}^{m-1}(-1)^i\alpha^{ki}\zeta(2ki+1,a)\beta^{k(m-i)}\zeta(2k(m-i)+1,b) \nonumber\\
		&-\pi^{2k}\beta^{k(m-1)}\sum_{n=0}^{\infty}\frac{1}{(n+b)^{2k(m-1)+1}}\sum_{i=0}^{\infty}\frac{1}{(i+a)(\beta^k(i+a)^{2k} +\alpha^k(n+b)^{2k})} \nonumber\\
		& -\pi^{2k}(-\alpha)^{k(m-1)}\sum_{n=0}^{\infty}\frac{1}{(n+a)^{2k(m-1)+1}}\sum_{i=0}^{\infty}\frac{1}{(i+b)(\beta^k(i+b)^{2k} +\alpha^k(n+a)^{2k})}\nonumber
	\end{align}
	
	which is equivalent to equation \eqref{eqtransform_2}.

	\subsection{First Proof of Theorem \ref{kernelrelation}}
	On account of the absolute convergence of $\zeta(s,a)$ for $\Re(s)>1$, $k>1$ we have
	\[\sum_{n=0}^{\infty} \frac{\Phi_{\alpha}(n+b,a;1)}{(n+b)^{2km}} = \frac{1}{2 i \pi}\int_{(c)} \frac{\zeta(1-s,a)\zeta(2km+s,b)}{2 \sin\left(\frac{\pi s}{2k}\right)}\left(\frac{\alpha}{\pi}\right)^{-s} ds\]

	We now evaluate this integral by shifting the line of integration and using Cauchy's Residue Theorem. Consider the contour determined by the line segments $[c-iT,c+iT],[c+iT,d+iT],[d+iT,d-iT],[d-iT,c-iT]$ where $d = -c-2km+1$. Inside this domian, the integrand has a simple pole at $-2km+1$ due to $\zeta(2km+s,b)$ and a pole of order two at $s=0$. It also has simple poles at the integers $-2ki$ where $i \in \{1,\ldots,m\}$ due to sine term in the denominator. The residues at these poles are 
	\begin{align*}
		& R_{0} = \frac{1}{\pi}\left(\zeta(2km,b)(\log\left(\frac{\alpha}{\pi}\right)+\gamma_{0}(a))-\frac{\partial}{\partial s}\zeta(2km+s,b)\big|_{s=0}\right) \\
		& R_{-2km+1}=\left(\frac{\alpha}{\pi}\right)^{2km-1}\frac{(-1)^m\zeta(2km,a)}{2\sin\left(\frac{\pi}{2k}\right)} \\
		& R_{-2ki}=(-1)^i \left(\frac{\alpha}{\pi}\right)^{2ki}\frac{\zeta(2ki+1,a)\zeta(2k(m-i),b)}{\pi}
	\end{align*}
	Thus by Cauchy's Residue theorem we have 
	\begin{align*}
		&\frac{1}{2i\pi}\left[\int_{c-iT}^{c+iT}+\int_{c+iT}^{d+iT}+\int_{d+iT}^{d-iT}+\int_{d-iT}^{c-iT}\right] \frac{\zeta(1-s,a)\zeta(s-2km-1,b)}{2\sin\left(\frac{\pi s}{2k}\right)}\left(\frac{\alpha}{\pi}\right)^{-s}\mathrm{d}s\\
		&= R_{0}+R_{-2km+1}+\sum_{i=1}^{m} R_{-2ki}
	\end{align*}
	From elementary bounds on Hurwitz zeta and \ref{sine}, it can be seen that as $T \to \infty$, the integrals along horizontal segments tend to zero. Under the change of variables $ s \to -s-2km+1$ we have 
	\begin{align*}
		&\int_{(d)} \frac{\zeta(1-s,a)\zeta(2km+s,b)}{2\sin\left(\frac{\pi s}{2k}\right)}\left(\frac{\alpha}{\pi}\right)^{-s}\mathrm{d}s \\
		& = \left(\frac{\alpha}{\pi}\right)^{2km-1}(-1)^m\int_{(c)} \frac{\zeta(2km+s,a)\zeta(1-s,b)}{2\cos\left(\frac{\pi (s+k-1)}{2k}\right)}\left(\frac{\beta}{\pi}\right)^{-s}\mathrm{d}s
	\end{align*}
	However,
	\begin{align*}
		\int_{(c)} \frac{\zeta(2km+s,a)\zeta(1-s,b)}{2\cos\left(\frac{\pi (s+k-1)}{2k}\right)}\left(\frac{\beta}{\pi}\right)^{-s}\mathrm{d}s = \sum_{n=0}^{\infty}\frac{\Psi_{\beta}(n+a,b;k)}{(n+a)^{2m}}
	\end{align*}
	which proves Theorem \ref{kernelrelation}.
	
	\subsection{Second Proof of Theorem \ref{kernelrelation}}
	Choose 
	\[a_n = \frac{1}{n+a} \,\, x_n = -\frac{(n+a)^{2k}}{\alpha^k} \,\, b_n = 1 \,\, \text{and} \,\, y_n = \frac{(n+b)^{2k}}{\beta^k}\]
	so the the corrosponding zeta functions are 
	\[\zeta_{x,a}(i)= (-1)^i\alpha^i\zeta(2ki+1,a) \,\, , \,\, \zeta_{y,b}(m-i)=\beta^{k(m-i)}\zeta(2k(m-i),b)\]
	and zeta generating functions are 
	\[\psi_{x,a}(z) = -\sum_{n=0}^{\infty}\frac{1}{n+a} \frac{\alpha^k z}{(n+a)^{2k}+\alpha^k z} \,\, , \,\, \psi_{y,b}(z) = \sum_{n=0}^{\infty} \frac{\beta^k z}{(n+b)^{2k}-\beta^k z}\]
	
	Using Theorem \ref{main_thm} we have 
	\begin{align*}
		& \sum_{i=1}^{m-1}(-1)^i\alpha^{ki}\zeta(2ki+1,a)\beta^{k(m-i)}\zeta(2k(m-i),b)\nonumber\\
		& =\pi^{2k}(-\alpha^k)^{m-1}\sum_{n=0}^{\infty}\frac{1}{(n+a)^{2k(m-1)+1}}\sum_{i=0}^{\infty}\frac{1}{\beta^k(n+a)^{2k}+\alpha^k(i+b)^{2k}}\\
		&-\pi^{2k}\beta^{k(m-1)}\sum_{n=0}^{\infty} \frac{1}{(n+b)^{2k(m-1)}}\sum_{i=0}^{\infty}\frac{1}{(i+a)(\beta^k(i+a)^{2k} +\alpha^k(n+b)^{2k})}\nonumber\\
	\end{align*}
	which can be seen equivalent to equation \eqref{kernelrelation}.

	\section{Concluding Remarks and future directions}
	
	Integral transforms are useful in number theory as can be witnessed from the
	results obtained in this paper.
	In a recent paper A. Dixit et al \cite{dixitsquared} found a Ramanujan Type identity for squares of zeta by considering the kernel 
	\[\Omega(x) = \int_{(c)} \frac{\zeta(1-s,a)^2}{2k\cos\left(\frac{\pi s}{2}\right)} x^{-s} ds,\]
	which they call Koshliakov kernel first introduced by N.S. Koshliakov in \cite{koshalikov} stated as follows
	\begin{align*}
		&\left(-\beta^{2}\right)^{-N}\left\{ \zeta^{2}\left(2N+1\right)\left(\gamma+\log\left(\frac{\beta}{\pi}\right)-\frac{\zeta'\left(2N+1\right)}{\zeta\left(2N+1\right)}\right)+\sum_{n=1}^{+\infty}\frac{\tau_{0}(n)}{n^{2N+1}}\,\Omega\left(\dfrac{\beta^2 n}{\pi^2}\right)\right\}
		\\ &  -\left(\alpha^{2}\right)^{-N}\left\{ \zeta^{2}\left(2N+1\right)\left(\gamma+\log\left(\frac{\alpha}{\pi}\right)-\frac{\zeta'\left(2N+1\right)}{\zeta\left(2N+1\right)}\right)+\sum_{n=1}^{+\infty}\frac{\tau_{0}(n)}{n^{2N+1}}\,\,\Omega\left(\dfrac{\alpha^2 n}{\pi^2}\right)\right\}
		\\ & = 2^{4N}\pi\sum_{j=0}^{N+1}\frac{\left(-1\right)^{j}\mathcal{B}_{2j}^{2}\mathcal{B}_{2N+2-2j}^{2}}{\left(\left(2j\right)!\right)^{2}\left(\left(2N+2-2j\right)!\right)^{2}}\left(\alpha^{2}\right)^{j}\left(\beta^{2}\right)^{N+1-j}.
	\end{align*}
	We have restated their transformation in terms of the Koshliakov kernel $\Omega(x)$, which has equivalent expressions \cite{dixitsquared}.
	
	By considering the generalized Koshaliakov kernel 
	\[\Omega(x,a;k)=\int_{(c)}\frac{\zeta(1-s,a)^2}{2k\cos\left(\frac{\pi(s+k-1)}{2k}\right)} x^{-s} ds\]
	one can obtain two parameter generalization of A. Dixits identity. 
	However $\zeta(1-s,a)^2$ does not have simple coefficients so we consider the special case $a=\frac{1}{n}$ where $n \in \NN$. Similar to the odd zeta kernel one can define a new kernel \[\Omega\left(x,a;k\right)=\int_{(c)}\frac{\zeta\left(1-s,a\right)^2}{2k\sin\left(\frac{\pi s}{2k}\right)} x^{-s} ds\]
	to find Generalization of Theorem \ref{Atul_Herglotz_1}. This will be explored in a future publication.
	
	\section{Acknowledgements}
	
	The author would like to thank Christophe Vignat and Tanay Wakhare for their guidance and support throughout the completion of this work and for taking time to read a rough draft of this manuscript.

\end{document}